\documentclass{article}

\usepackage{amsthm}
\usepackage{amssymb}
\usepackage{amsmath}

\newtheorem{theorem}{Theorem}[section]
\newtheorem{lemma}[theorem]{Lemma}

\theoremstyle{definition}

\theoremstyle{remark}

\newtheorem{proposition}[theorem]{Proposition}

\numberwithin{equation}{section}

\usepackage{natbib}
\usepackage{color}
\usepackage{chngcntr}
\usepackage{soul}

\def\bK{{\mathbb K}}

\def\bZ{{\mathbb Z}}

\def\bN{{\mathbb N}}
\def\bF{{\mathbb F}}

\def\bR{{\mathbb R}}

\def\bS{{\mathbb S}}

\def\bK{{\mathbb K}}

\def\bZ{{\mathbb Z}}

\def\bN{{\mathbb N}}
\def\bF{{\mathbb F}}

\def\bR{{\mathbb R}}

\def\bS{{\mathbb S}}

\def\cH{\mathcal{H}}

\def\R{\mathbb{R}}

\def\AA{B}

\allowbreak

\newcommand{\balpha}{ \mbox{\boldmath $ \alpha $} }

\newcommand{\bmu}{ \mbox{\boldmath $\mu$} }

\newcommand{\bsigma}{ \mbox{\boldmath $\sigma$} }

\renewcommand{\S}{\mathbb{S}}

\newcommand{\bB}{\textbf{B}}
\newcommand{\bk}{\textbf{k}}
\newcommand{\bh}{\textbf{h}}
\newcommand{\bu}{\textbf{u}}

\newcommand{\bx}{\textbf{x}}

\def\cH{\mathcal{H}}

\def\R{\mathbb{R}}

\def\AA{B}

\begin{document}

\title{Multivariate Gaussian Random Fields over Generalized Product Spaces involving the Hypertorus}


\maketitle

\begin{center}
	Fran\c{c}ois Bachoc \\
	Department of Mathematics, Universit{\'e} Paul Sabatier, Toulouse, France \\
	Ana Paula Peron\footnote{The second author was partially supported by FAPESP} \\
	Department of Mathematics, ICMC,  University of S\~ao Paulo, S\~ao Carlos, Brazil \\
	Emilio Porcu\footnote{The third author was supported in part by ADD} \\
	Department of Mathematics, Khalifa University, The United Arab Emirates, $\&$ School of Computer Science and Statistics, Trinity College Dublin
\end{center}

\begin{abstract}
The paper deals with multivariate Gaussian random fields defined over generalized product spaces that involve the hypertorus. The assumption of Gaussianity implies the finite dimensional distributions to be completely specified by the covariance functions, being in this case matrix valued mappings. \\
We start by considering the spectral representations that in turn allow for a characterization of such covariance functions. We then provide some methods for the construction of these matrix valued mappings. Finally, we consider strategies to evade radial symmetry (called isotropy in spatial statistics) and provide representation theorems for such a more general case.  
~ \\
{This paper is dedicated to Professor Yadrenko, who has largely inspired our research since the times of our PhD studies.} 
\end{abstract}

\section{Introduction}\label{Introduction}
\setcounter{equation}{0}

{\vspace{0.3cm}}

\def\S{\mathbb{S}}
\def\T{\mathbb{T}}
\def\k{d_1^{\prime}}
\def\x{x^{\prime}}

\subsection{Context}

This paper considers multivariate Gaussian random fields defined continuously over product spaces involving the torus, or generalizations of the torus obtained through the product of hyperspheres cross $d$-dimensional Euclidean spaces. Several applications from applied sciences motivate our work on these geometries, and we discuss some of these settings throughout the paper. \\
More specifically, we consider generalized product spaces, that we define here as the Cartesian product of the $d$-dimensional Euclidean space, $\R^d$, with the hypertorus of $(d_1,d_2)$ dimensions, being in turn the product of two unit hyperspheres that are embedded in $\R^{d_1+1}$ and $\R^{d_2+1}$, respectively. In particular, we consider multivariate Gaussian fields (vector valued fields) as in typical applications, see for instance  the survey in \cite{PAF2016}, several variables are observed at the same space-time sampling points. Hence, there is a need to take into account not only the space-time dependence of one variable, but also the cross dependence between any pair of random fields that are observed over the space-time domain. 

Most of the literature in the last $30$ years has considered {\em space} as a planar domain - the Euclidean space $\R^d$ - while {\em time} has been typically considered as linear \citep{porcu30}.
The main reason for considering the above generalized spaces based on the hypertorus as {\em space-time} index sets is that these allow for more generality for a variety of situations. For univariate (scalar valued) random fields, some special cases of these generalized spaces are provided by \cite{PW1}. For instance, continuous-time data often exhibit multiple sources of seasonality ({\em e.g.}, daily and weekly). One potential strategy to address these seasonalities is by incorporating them into a covariance structure by wrapping time over the product of two circles, $\S^1 \times \S^1$, see \cite{shirota} and \cite{white2019a} on continuous-time monitoring of ground-level ozone concentrations. 

As a second possible application, one might consider data over a large section of the Earth where the data exhibit strong seasonality (perhaps on a daily, weekly, or annual scale) or are indexed by a direction (\textit{e.g.}, a bird flying southward). The planet is often approximated by a sphere $\S^2$ with a given radius \citep{PAF2016}, and, as discussed, time can be wrapped in a circle to capture seasonality over time. Alternatively, the direction of the data could be indexed by an angle on the circle. Thus, in either case, a Gaussian random field with an index set of $\S^2 \times \S^1$ would be a natural choice to represent phenomena in \textsc{Space}$\times$\textsc{Time} or \textsc{Space}$\times$\textsc{Direction} \citep{mastrantonio2}.

Gaussian random fields have finite dimensional distributions that are completely specified by their mean, and their covariance function. Covariance functions are positive definite, and ensuring such a requirement for a given candidate mapping is a nontrivial task. For this paper, such a task is overly complicated by the fact that (a) for vector valued random fields, the covariance function is a matrix valued function, and the mathematical machinery behind {\em positivity} for the matrix valued case is very challenging \citep[see, for instance][]{cramer1940theory}. Also, (b) generalized spaces have been considered to a very limited extent - see the subsequent literature discussion. 

\subsection{Literature}

The last ten years have seen an unprecedented interest in random fields that are defined over non planar surfaces. For a thorough account, the reader is referred to \cite{porcu30}, with the references therein.

Recently, there has been some interest in modeling random fields defined over $\R^2 \times \S^1$, where $\R^2$ is {\em space} and $\S^1$ is time wrapped over the circle. For instance, \cite{shirota} challenge on daily crime events in San Francisco (space, a compact set of $\R^2$) in $2012$ (time, wrapped over the circle). The same modus operandi is adopted by \cite{mastrantonio}, who consider Bayesian hierarchical modeling where seasonality is modeled through conditioning sets \citep[see also][]{white2019a}.

Let $d,d_1$ and $d_2$ be positive integers. 
Random fields defined over the product space $\S^{d_1} \times \S^{d_2}$ have only been considered to a limited extent. In mathematical analysis, \cite{mene2, mene4, mene3} and \cite{mene1} provide conditions for positive definiteness. \cite{berg2016schoenberg} instead considered the space $\S^d \times G$, for $G$ being a locally compact group. \cite{PBG16} and \cite{white} considered Gaussian fields defined over $\S^d$ cross linear time. The reader is referred to \cite{PAF2016} for a recent review on the subject. Recent generalizations on product spaces involving hyperspheres have been proposed in \cite{estrade}, and, in more abstract settings, recent contributions can be found in \cite{berg1} and \cite{berg2}.

\def\AA{\textsc}


\subsection{The Problems, and our Contribution}
We start by highlighting some aspects where the literature currently lacks. \\
{\bf A.} Multivariate Gaussian random fields defined over generalized spaces are unexplored. \\
{\bf B.} Characterizations theorems - that is, spectral representations - for the covariance functions associated with multivariate Gaussian random fields are not available in the current literature. \\
This paper is devoted to some challenges related to this topic. \\ 
{\bf 1.} We consider matrix valued covariance functions that are componentwise isotropic over generalized spaces. That is, they are radially symmetric in $\R^d$, and they depend on the dot product of each sphere involving the hypertorus. For this setting, \\
{\bf 1.a) } we provide a spectral representation theorem that ensures necessary and sufficient conditions for a given continuous candidate matrix valued function to be positive definite; \\
{\bf 1.b)} we provide constructive criteria that allow to build new parametric classes of matrix valued covariance functions under this setting{, in particular of the Mat\'ern type}.
\\
{\bf 2.} We relax the hypothesis of radial symmetry in $\R^d$ and provide spectral representations and convergence conditions for the spectral expansions to be well defined. \\ 

The plan of the paper is the following: {Section \ref{sec2} provides relevant background material. Section \ref{sec:theoretical:results} 
		 first provides a spectral representation of covariance functions on generalized spaces and then provides new general classes of covariance functions, in particular new classes of the Mat\'ern type. Finally, Section \ref{sec:4} studies a class of matrix covariance functions on generalized spaces that are not radially symmetric in $\R^d$, providing a characterization.}

\section{Background Material}
\label{sec2}

\subsection{Background from Harmonic Analysis}

Let $\mathbb{S}^n$ be the unit sphere of dimension $n \geq1$, embedded in $\bR^{n+1}$, and let $\Delta$ be the Laplace-Beltrami operator defined over $\mathbb{S}^n$ {(see Section 2.2 of \cite{MR1641900})}.
 The spectrum of $\Delta$ is discrete, real and non-positive, with eigenvalues given by $\lambda_k=-k(k+n-1)$ for $k\geq 0$ \citep{szego}. Let $\mathcal{H}_k$ be the eigenspace that corresponds to the eigenvalue $\lambda_k$. It is known that \citep{MR1641900} 
$$L^2(\bS^n)=\bigoplus_{k=0}^\infty \cH_k.$$
Let 
$\{Y_{k,j}:\;k\in\mathbb{N}_0,\;j=1,2,\dots,\kappa_{n,k}\}$ be an orthonormal basis of $\mathcal{H}_k$ having uniquely determined dimension, $\kappa_{n,k}$, defined as
\begin{equation}\label{apr1}
	\kappa_{n,k}=\text{dim} \mathcal{H}_k= (2k+n-1)\frac{(k+n-2)!}{k!(n-1)!},
\end{equation}
where this expression has to be understood as bein identically equal to $1$ when $k=0$ and $n=1$. 
From (\ref{apr1}) it becomes apparent that there exists a constant $c_n\ge1$ such that
\begin{equation}\label{dimHksim}
	c_n^{-1} (k+1)^{n-1} \leq \kappa_{n,k}  \leq c_n (k+1)^{n-1}.
\end{equation}
The functions $Y_{k,j}$ are the well-known spherical harmonics on $\mathbb{S}^n$,  and they satisfy the addition formula \citep{marinucci-peccati, yadrenko1983spectral}
\begin{equation}\label{apr2}
	\sum_{j=1}^{\kappa_{n,k} }Y_{k,j}(x)\overline{Y_{k,j}(y)}=\frac{\kappa_{n,k}}{\omega_n} C_k^{(n-1)/2}(\langle x,y\rangle), \qquad x,y \in \S^n,
\end{equation}
where $\omega_n={(2\pi^{(n+1)/2})}/{\Gamma((n+1)/2)}$ is the total area of $\mathbb{S}^n$, with 
$$C_k^{n}(r)= \frac{c_{k}^{\frac{n-1}{2}}(r)}{c_{k}^{\frac{n-1}{2}}(1)},  \qquad r \in [-1,1],$$ where $c_k^{\lambda }$ is a Gegenbauer polynomial of degree $k$ and order $\lambda > -1/2$ \citep{szego}. Fundamentally, if $\lambda>-1/2$, then there exists a constant $c=c_{\lambda}>0$ such that
\begin{equation}
	\label{Gc}
\sup_{r \in [-1,1]}	c_{k}^{\lambda}(r)\le c(k+1)^{2\lambda-1},\quad\text{for every}\;\;k=0,1,\ldots.
\end{equation}

\subsection{Gaussian Fields on Generalized Product Spaces involving the Hypertorus}

Let $d_1,d_2$ be positive integers. We define the $(d_1,d_2)$-hypertorus with index $(d_1,d_2)$ through the identity 
$$ \T^{d_1,d_2} := \S^{d_1} \times \S^{d_2} = \{ \mathbf{x}=(x_1,x_2) \; : x_i \in \R^{d_i+1}, \; \|x_i\|_{d_i+1}=1,\;i=1,2\}, $$
where $\|\cdot\|_{d_i+1}$ is the Euclidean norm on $\bR^{d_i+1}$. 
The name {\em hypertorus} is due to $\T^{1,1}= \S^1 \times \S^1$, the product of two circles, being isomorphic to the classical torus. 

\def\bZ{\boldsymbol{Z}} 
\def\bu{\mathbf{u}}
\def\bK{\boldsymbol{K}}
\def\bh{\boldsymbol{h}}
\def\B{\boldsymbol{B}}
\def\buu{\mathbf{u}^{\prime}}
\def\bmu{\boldsymbol{\mu}}
\def\balpha{\boldsymbol{\alpha}}
\def\bsigma{\boldsymbol{\sigma}}

Let $p,d$ be positive integers, and let $\top$ be the transpose of a column vector. The paper considers zero mean $p$-variate (or simply, multivariate) Gaussian random fields $\{\bZ(\mathbf{x},\mathbf{u})= (Z_1(\mathbf{x},\bu),\ldots,Z_p(\mathbf{x},\bu))^{\top}, \; \mathbf{x} \in \T^{d_1,d_2}, \bu \in \R^d \}$ and focuses on their covariance function $\bK_{\boldsymbol{Z}}: \left ( \T^{d_1,d_2} \times \R^d \right ) \times \left ( \T^{d_1,d_2} \times \R^d \right ) \to \R^{p \times p}$, having elements $K_{ij,\boldsymbol{Z}}$, $i,j=1,\ldots, p$, that are defined as 

$$ K_{ij,\boldsymbol{Z}} \big ( (\mathbf{x}, \bu),(\mathbf{x}^{\prime},\bu^{\prime}) \big ) = \mathbb{E} \Big ( Z_{i}(\mathbf{x},\bu) Z_j(\mathbf{x}^{\prime},\bu^{\prime}) \Big ), \qquad \mathbf{x}, \mathbf{x}^{\prime} \in \T^{d_1,d_2}, \bu,\bu^{\prime} \in \R^d. $$

For spectral representations for Gaussian fields, the reader is referred to \cite{malyarenko2018tensor}.
Covariance functions  are a linear measure of dependence between any pair of random variables $Z_i(\mathbf{x},\bu)$, $Z_j(\mathbf{x}^{\prime},\bu^{\prime})$ located at points $(\mathbf{x},\bu),(\mathbf{x}^{\prime},\bu^{\prime}) \in \T^{d_1,d_2} \times \R^d$. The function $\bK_{\boldsymbol{Z}}$ is positive definite \citep{schoen42}: for any collection $\{ \boldsymbol{a}_k\}_{k=1}^N \subset \R^{p N}$ and points $\{(\mathbf{x}_k,\bu_k)\}_{k=1}^N \in {(} \T^{d_1,d_2} \times \R^d {)^N}$, 
$$ \sum_{k,h=1}^{N} \boldsymbol{a}_k^{\top} \bK_{\boldsymbol{Z}} \Big ((\mathbf{x}_k,\bu_k),(\mathbf{x}_h,\bu_h) \Big) \boldsymbol{a}_h \ge 0. $$
Building covariance functions is mathematically challenging and simplifying assumptions are often required for modeling, estimation, and prediction. Section \ref{sec:theoretical:results} assumes that there exists a continuous (continuity is intended as pointwise) mapping $\bK : [-1,1]^2 \times [0,\infty) \to \R^{p \times p}$ with $K_{ii}(1,1,0)=1$,
such that 
\begin{equation}
	\label{menegatto}  \bK_{\boldsymbol{Z}} \Big ( (\mathbf{x},\bu), (\mathbf{x}^{\prime},\bu^{\prime}) \Big)  = \Sigma^{1/2}   \bK( s,r,h  ) \Sigma^{1/2}, \qquad 
\end{equation} 
with $s := \langle x_1,{x_1^\prime} \rangle_1,\,  r :=  \langle x_2,{x_2^\prime} \rangle_2 $,\,  $h:= \| { \bu - \bu^\prime }\|$
for {$\mathbf{x}=(x_1,x_2), \mathbf{x}^\prime=(x_1^\prime,x_2^\prime)$,} $x_i,x_i^{\prime} \in \S^{d_i}$ and where $\langle \cdot, \cdot \rangle_i$ denotes the classical dot product on $\R^{d_i+1}$ $i=1,2$. Here, $\Sigma^{1/2}$ is a diagonal matrix such that $(\Sigma^{1/2})^\top \Sigma^{1/2} =\Sigma= {\rm diag} (\sigma_1, \ldots,\sigma_p)$ with $\sigma_i$ denoting the variance of $Z_i$.

We call the covariance functions $\bK_{\boldsymbol{Z}}$ componentwise isotropic. In Equation (\ref{menegatto}), $s$ and $r$ are the cosine of the great circle distances taken over $\S^{d_1}$ and $\S^{d_2}$, respectively. How to relax the hypothesis of componentwise isotropy will be discussed in Section \ref{sec:4}. Throughout the rest of the paper, we say that $\bK$ is positive definite to mean that the composition of $\bK$ with the triplet $(\langle \cdot,\cdot \rangle_1,  \langle \cdot,\cdot \rangle_2, \|\cdot \|)$ is a positive definite matrix valued function. This abuse of notation will allow to have a simplified exposition of the results coming subsequently.

\subsection{Conditionally Negative Definite Functions}

Let $X$ be any nonempty set.
A matrix valued function $\boldsymbol{\gamma} : X \times X \to \R^{p \times p}$ is a cross variogram if it is conditionally negative definite, that is 
$$ \sum_{k,h=1}^N \boldsymbol{a}_k^{\top} \boldsymbol{\gamma}(\boldsymbol{x}_k,\boldsymbol{x}_h) \boldsymbol{a}_h \le 0, $$
for all $\{ \boldsymbol{x}_k \}_{k=1}^N$ being points in $X$ and {$\{\boldsymbol{a}_k\}_{k=1}^N$} being a collection of vectors summing up to the zero vector. A similar definition applies to matrices $\boldsymbol{A}= [a_{ij}]_{i,j=1}^p$. 
We observe that $\gamma$ is conditionally negative definite if and only if $-\gamma$ is conditionally positive definite (that is, the sign in the inequality above is changed by greater or equal): we refer \cite{2004P-valdir-cpd-eucl, 2006P-valdir-claude-cpd-dotkernel} for results related to characterizations  of conditional positive definiteness when $p=1$ and $X$ is respectively $\S^d$ or $\mathbb R^d$,  and \cite{guella-men-cpd2020} when $X=\R^d$, $p\geq1$ and $N$ is fixed.



\section{Theoretical Results} \label{sec:theoretical:results} 

\subsection{Spectral Representations}

We start with a characterization theorem that will provide the basis to subsequent findings. The result provides a spectral representation for all functions $\bK$ that are defined through Equation (\ref{menegatto}).
\begin{theorem} \label{mauro_avoid!}
	Let $\bK_{\boldsymbol{Z}}$ be the mapping defined through Equation (\ref{menegatto}). Then, $\bK_{\boldsymbol{Z}}$ is positive definite if and only the mapping $\bK$ on the right hand side of (\ref{menegatto}) is uniquely determined through the expansion
	\begin{equation}\label{CZ}
		\bK(s,r,h)=\sum_{\boldsymbol{k}=(k_1,k_2)\in\bN_0^2} \B^{(d_1,d_2)}_{\boldsymbol{k}}(h) C_{k_1}^{\frac{d_1-1}{2}}(s) C_{k_2}^{\frac{d_2-1}{2}}(r), \qquad s,r \in [-1,1], \; h \ge 0,
	\end{equation}
	where $\{ \B^{(d_1,d_2)} _{\boldsymbol{k}}(\cdot)\}$ is a sequence of continuous functions from $[0,\infty)$ into $\R^{p \times p}$ such that 
	\begin{enumerate}
		\item $\sum_{\boldsymbol{k}} \B_{\boldsymbol{k}}^{(d_1,d_2)} (0) = \B$, where $\B$ is a matrix with ones in the diagonal. 
		\item For every $\boldsymbol{k} \in \mathbb{N}_0^2$, the function $\R^d \ni \bh \mapsto \B_{\boldsymbol{k},Z}(\bh)$, defined as 
		$$ \B_{\boldsymbol{k},Z}(\bh) = \B_{\boldsymbol{k}}^{(d_1,d_2)}(\|\bh\|), \qquad \bh \in \R^d, $$
		is the covariance function of a multivariate Gaussian random field $\boldsymbol{Z}_{\boldsymbol{k}}$. 
	\end{enumerate}
\end{theorem}
The proof is based on the following 
\begin{lemma} \label{Lula_free_again}
	Let $\bK: [-1,1] \times [-1,1] \times [0,\infty) \to \R^{p \times p}$. Then, $\bK$ is positive definite if and only if 
	\begin{equation} \label{eq:lemma1}
		\int_{\T^{d_{1},d_2}} \int_{\R^d}  \int_{\T^{d_{1},d_2}} \int_{\R^d}  {\rm d} \bmu \big ( x_1,x_2, \bu \big )^{\top } \bK \Big ( \langle x_1,x_1^{\prime}  \rangle, \langle x_2, x_2^{\prime} \rangle, \|\bu- \buu \|    \Big )  \overline{{\rm d}\bmu \big (x_1^{\prime},x_2^{\prime}, {\bu^\prime} \big )} \ge 0, 
	\end{equation} 
	for every complex-valued Radon measure $\bmu :\T^{d_1,d_2} \times \R^d \to \mathbb{C}^{p}$.
\end{lemma}
{\bf Proof.}
It is clear that if $\bK$ satisfies \eqref{eq:lemma1}, then $\bK$ is positive definite.
Let $\bK$ be positive definite. Then, $\bK$ is bounded \citep{berg2016schoenberg}. Take $\bmu$ as 
$$ \bmu(x_1,x_2,\bu) = \sum_{j=1}^n \balpha_{j} \delta_{(x_{1,j},x_{2,j},\bu_j)}(x_1,x_2,\bu), 
 \qquad (x_{1,j},x_{2,j},\bu_j) \in \T^{d_{1},d_2} \times \R^d, $$
for $\balpha_j \in \R^p$  and with $\delta$ denoting the Kronecker delta. Then, the left hand side of Equation (\ref{eq:lemma1}) can be written as 
$$ \sum_{k,\ell}^n \balpha_k^{\top} \bK \Big ( \langle x_{1,k},x_{1,\ell} \rangle, \langle x_{2,k}, x_{2,\ell} \rangle, \|\bu_k- \bu_{\ell} \|    \Big ) \balpha_{\ell}, $$
which is apparently nonnegative because $\bK$ is positive definite. Now, for a general $\bmu$, consider the number 
$$ I:= \int_{\T^{d_{1},d_2}} \int_{\R^d}  \int_{\T^{d_{1},d_2}} \int_{\R^d}  {\rm d} \bmu \big (x_1,x_2, \bu \big )^{\top } \bK \Big ( \langle x_1,x_1^{\prime}  \rangle, \langle x_2, x_2^{\prime} \rangle, \|\bu- \buu \|    \Big )  \overline{{\rm d} \bmu \big (x_1^{\prime},x_2^{\prime}, \bu \big )}. $$
We will prove that $I\ge 0$ by showing that for every $\epsilon >0$ there exists $J \ge 0$ such that $\Big | I-J\Big | < \epsilon || \bmu ||^{2}$, with $|| \cdot ||$ denoting the total variation of a complex valued measure with image in $\mathbb{C}^p$ \citep{rudin}. The rest of the proof follows straight by invoking the same arguments as in the proof of Lemma 4.3 of \cite{berg2016schoenberg}. 
\hfill $\Box$ \\
\\
{\bf Proof of Theorem \ref{mauro_avoid!}. } We start by supposing that $\bK$ is positive definite. Then, we apply Lemma \ref{Lula_free_again} under the choice 
$$ \bmu = \sigma_{d_1} \otimes \sigma_{d_2} \otimes \bsigma, $$
where $\otimes $ denotes tensor product, $\sigma_n$ is the Lebesgue measure for $\S^n$, and $\bsigma$ is an arbitrary Radon vector valued measure with compact support. Thus, we have that the integral
\begin{eqnarray}
	&& \int_{(\T^{d_{1},d_2})^2} \int_{\R^{2d}} {\rm d}  \sigma_{d_1}(x_1) {\rm d}  \sigma_{d_2}( x_2) {\rm d}  \bsigma (\bu)^{\top} \bK \Big ( \langle x_1,x_1^{\prime}  \rangle, \langle x_2, x_2^{\prime} \rangle, \|\bu- \buu \|    \Big )  \nonumber \\ 
	&& \overline{{\rm d} \sigma_{d_1}( x_1^{\prime})} \overline{{\rm d} \sigma_{d_2}( x_2^{\prime}) } \overline{{\rm d} \bsigma( \bu^{\prime})}
	  \label{berg1}
\end{eqnarray}   is nonnegative. 
Using the tensor product measure as above allows to write 
\begin{eqnarray*} && \int_{\T^{d_{1},d_2}}  \int_{\T^{d_{1},d_2}}  {\rm d} \sigma_{d_1}(x_1)  {\rm d}  \sigma_{d_2}( x_2) \bK \Big ( \langle x_1,x_1^{\prime}  \rangle, \langle x_2, x_2^{\prime} \rangle, \|\bu- \buu \|    \Big ) \overline{{\rm d} \sigma_{d_1}( x_1^{\prime})} \overline{{\rm d} \sigma_{d_2}( x_2^{\prime}) } \\
	 &=& \qquad \omega_{d_1} \omega_{d_2} \int_{\T^{d_{1},d_2}} \bK\Big ( \langle {\rm e}_1, \eta_1  \rangle, \langle {\rm e}_1, \eta_2 \rangle, \|\bu- \buu \|    \Big ) {\rm d} \sigma_{d_1}(\eta_1) {\rm d}  \sigma_{d_2}( \eta_2), \end{eqnarray*}  with ${\rm e}_1=(1,0,\ldots,0)^{\top}$ Hence, (\ref{berg1}) can be rewritten as 
\begin{equation}
	\label{berg2} 
	\int_{\T^{d_{1},d_2}} \int_{\R^d}  \int_{\R^d}  {\rm d}  \bsigma ( \bu)^{\top} \bK \Big ( \langle {\rm e}_1,\eta_1  \rangle, \langle {\rm e}_1, \eta_2\rangle, \|\bu- \buu \|    \Big )  \overline{{\rm d}  \bsigma(\bu^{\prime})} {\rm d}  \sigma_{d_1}( \eta_1) {\rm d}  \sigma_{d_2}( \eta_2)\ge 0. 
\end{equation}
We now apply (\ref{berg2}) to the matrix valued function $\bK(\cdot,\cdot,\cdot)  C_{k_1}^{\frac{d_1-1}{2}}(\cdot) C_{k_2}^{\frac{d_2-1}{2}}(\cdot)$
(that is positive definite because $\bK$ is positive definite by hypothesis, and the Gegenbauer polynomials are positive definite by construction) to obtain
\begin{eqnarray}
	\label{berg22} 
	&& \int_{\T^{d_{1},d_2}} \int_{\R^d}  \int_{\R^d}  {\rm d}  \bsigma (\bu)^{\top} \bK \Big ( \langle {\rm e}_1,\eta_1  \rangle, \langle {\rm e}_1, \eta_2\rangle, \|\bu- \buu \|    \Big ) \nonumber \\ &\times & \qquad  C_{k_1}^{\frac{d_1-1}{2}}(\langle {\rm e}_1 , \eta_1 \rangle ) C_{k_2}^{\frac{d_2-1}{2}}(\langle {\rm e}_1 , \eta_2 \rangle) \overline{{\rm d} \bsigma( \bu^{\prime})} {\rm d} \sigma_{d_1}(\eta_1) {\rm d}  \sigma_{d_2}( \eta_2)\ge 0. 
\end{eqnarray}
We can now define the sequence of fuctions $\B_{\boldsymbol{k}}^{(d_1,d_2)}$ through the identity
\begin{eqnarray*}
	\B_{\boldsymbol{k}}^{(d_1,d_2)}(h) &=&  \frac{\kappa_{d_1,k_1} \kappa_{d_2,k_2}}{\omega_{d_1} \omega_{d_2}} \int_{\T^{d_1,d_2}} \bK \Big ( \langle {\rm e}_1,\eta_1  \rangle, \langle {\rm e}_1, \eta_2\rangle, {h}    \Big ) \\
	&\times & \qquad C_{k_1}^{\frac{d_1-1}{2}}(\langle {\rm e}_1 , \eta_1 \rangle ) C_{k_2}^{\frac{d_2-1}{2}}(\langle {\rm e}_1 , \eta_2 \rangle) {\rm d}  \sigma_{d_1}(\eta_1){\rm d}  \sigma_{d_2}( \eta_2),
\end{eqnarray*}
$h \ge 0$. Equation (\ref{berg22}) shows that $\B_{\boldsymbol{k}}^{(d_1,d_2)}$ is positive definite on $\R^d$ for every fixed $\bk=(k_1,k_2)\in\bN_0^2$. 
Clearly, the identity above can be simplified into 
\begin{eqnarray*}
	&& \B_{\boldsymbol{k}}^{(d_1,d_2)}(h)  \\ &=&\frac{\kappa_{d_1,k_1} \kappa_{d_2,k_2}}{\omega_{d_1} \omega_{d_2}} \int_{-1}^{1} \int_{-1}^{1} \bK \Big ( r, s, h    \Big ) C_{k_1}^{\frac{d_1-1}{2}}(r) C_{k_2}^{\frac{d_2-1}{2}}(s) 
\frac{(1-r^2)^{d_1/2-1}}{ \mathcal{B}( d_1/2 , 1/2 )}
\frac{(1-s^2)^{d_2/2-1}}{ \mathcal{B}( d_2/2 , 1/2 )}
	  {\rm d} r {\rm d} s.
  \end{eqnarray*}
Here, $\mathcal{B}$ denotes the Beta function. 
The rest of the proof follows the same arguments as in \cite{berg2016schoenberg} and is thus omitted. \hfill $\Box$

\subsection{New Classes of Covariance Functions} \label{sec3}

The previous theorem is actually the key to prove the following results. 

\begin{theorem} \label{thm:bolsonaro1}
	Let $\boldsymbol{\gamma}: [-1,1] \times [0,\infty) \to \R^{p \times p}$ such that
	$- \boldsymbol{\gamma}(s,h)$ is a cross variogram over $\left ( \S^{d_2} \times \R^d \right )^2$. Additionally, suppose the elements $\gamma_{ij}$ of $\boldsymbol{\gamma}$ are all {strictly positive valued}. Let $\bK$ have elements $K_{ij}$ that are defined as
	\begin{equation}
		\label{Bolsonaro1} 
		K_{ij}(s,r,h) = \frac{1}{ \gamma_{ij}(r,h)} + \frac{\pi}{2} \frac{\sinh  \sqrt{\gamma_{ij}(r,h) (\pi - \arccos s)}}{\sinh \sqrt{\gamma_{ij}(r,h)}}, \qquad s,r \in [-1,1], \quad h \ge 0.
	\end{equation}
	Then, $\bK(s,r,h)$ is positive definite on $\T^{d_1,d_2} \times \R^d$ for $d_1=1$ and for all ${d,}d_2=1,2,\ldots$.
\end{theorem}
{\bf Proof.} We start by noting that the Online Supplement in \cite{porcu2019axially} proves the identity
$$   \frac{1}{ \gamma_{ij}(r,h)} + \frac{\pi}{2} \frac{\sinh  \sqrt{\gamma_{ij}(r,h) (\pi - \arccos s)}}{\sinh \sqrt{\gamma_{ij}(r,h)}} = \sum_{k=0}^{\infty} \frac{\cos \Big ( k \arccos s \Big ) }{\Big (k^2 + \gamma_{ij}(r,h) \Big ) }, \qquad s,r \in [-1,1], \quad h \ge 0. $$

Hence, the function $\bK$ as in (\ref{Bolsonaro1}) can be written as 
$$ \bK(s,r,h) = \sum_{k}\cos \Big ( k \arccos s \Big ) \boldsymbol{B}_k(r,h), \qquad s,r \in [-1,1], h \ge 0, $$
with $\{ \bB_k \}$ being a sequence of matrix valued functions having entries $B_{ij,k}$ defined as 
$$B_{ij,k}(r,h) =  \frac{1}{\Big (k^2 + \gamma_{ij}(r,h) \Big ) }, \qquad r \in [-1,1], \quad h \ge 0. $$
Hence, the proof is completed if we prove that $\{  \bB_k(\cdot,\cdot) \}_k$ is a sequence  of positive definite functions that are summable at $(r,h)=(1,0)$. Indeed, we have that, for every $k=0,1,\ldots$,
$$ B_{ij,k}(r,h) = \int_{0}^{\infty} {\rm e}^{-\xi k^2} {\rm e}^{- \xi \gamma_{ij}(r,h)} {\rm d} \xi. $$
Since $- \boldsymbol{\gamma}$ is a cross variogram, a direct application of Theorem 1 in \cite{schlather2010some} shows that $\bB_k$ is positive definite for every fixed $k$. Summability of $\bB_k(1,0)$ si obtained through direct inspection. The proof is completed. \hfill $\Box$

We now define the parametric class of functions $\mathcal{F}(\cdot; \alpha,\tau,\nu)$ through the identity 
\begin{equation} 
	\label{F-class} \mathcal{F}(t; \alpha, \tau, \nu) = \frac{\mathcal{B}(\alpha,\nu+\tau)}{\mathcal{B}(\alpha,\nu)} \;_{2}F_{1} \left (  \tau,\alpha; \alpha+\nu+\tau; -t \right ), \qquad t \in [-1,1], 
\end{equation}
where $\;_{2}F_{1}$ is an hypergeometric function \citep{olver2010nist}, and where $\alpha,\nu,\tau$ are strictly positive. Here, $\mathcal{B}$ denotes the Beta function. 

\begin{proposition}
	Let $\boldsymbol{C}: [-1,1]\times [0,\infty) \to \R^{p \times p }$ be positive definite in $\S^{d_2} \times \R^d$ with entries $C_{ij}$ such that $|C_{ij}(\cdot)| \le 1$. Let $[\alpha_{ij}]$, $[\beta_{ij}]$ be conditionally negative definite matrices.  
	Then, 
	the function $\bK: [-1,1] \times [-1,1] \times [0,\infty) \to \R^{p \times p}$ having entries $K_{ij}$ defined as 
	$$ K_{ij}(s,r,h) = \mathcal{B}(\alpha_{ij},\beta_{ij}) \mathcal{F}( s C_{ij}(r,h); \alpha_{ij}, \tau, \nu_{ij}), $$
	is positive definite on $\T^{d_1,d_2} \times \R^d$ for all $d_1=1,2, \ldots$. 
\end{proposition}

{\bf Proof.} We start by considering the class 
\begin{equation} \label{neg-bin}
	\mathcal{N}_{\delta,\tau}(s) = \left ( \frac{1-\delta}{1-\delta s} \right )^{\tau}, \qquad s \in [-1,1], \quad \delta \in (0,1), \; \tau>0.
\end{equation}
Direct inspection shows that $\mathcal{N}_{\delta,\tau}$ admits the followig expansion:
$$ \mathcal{N}_{\delta,\tau}(s)  = \sum_{k=0}^{\infty} s^k b_{k}(\delta,\tau), \qquad s \in [-1,1], $$
where $b_{k}(\delta,\tau) = \binom{k+\tau-1}{k} \delta^{k} (1-\delta)^{\tau}$. Hence,
 the function 
$$ \mathcal{N}_{i,j;\delta,\tau}(s)  = \sum_{k=0}^{\infty} s^k b_{k}(\delta,\tau)C_{ij}(r,h)^k, \qquad s,r \in [-1,1], h \ge 0, $$
is the $(i,j)$-th entry of a positive definite function, $\boldsymbol{\mathcal{N}}$. We now apply the arguments in theorem 1 of \cite{alegria2021mathcalffamily} to infer that 
$$ \mathcal{F}( s C_{ij}(r,h); \alpha_{ij}, \tau, \nu_{ij}) = \int_{0}^1 \mathcal{N}_{ij;\delta,\tau}(s) \frac{\delta^{\alpha_{ij}-1} (1-\delta)^{\beta_{ij}-1}}{\mathcal{B}(\alpha_{ij},\beta_{ij})} {\rm d} \delta. $$
We now observe that, since $[\alpha_{ij} ]_{i,j=1}^p $ and $[\beta_{ij}]_{i,j=1}^p$ are conditionally negative definite matrices, we get $\delta^{\alpha_{ij}} = \exp (\alpha_{ij} \log \delta)$, proving that the matrix $[\delta^{\alpha_{ij}}]_{i,j=1}^p$ is positive definite{, from a theorem by Schoenberg (see for instance Theorem A.3 in \cite{bachoc2017gaussian} or \cite{berg1984harmonic})} . So is the matrix $[(1-\delta)^{\beta_{ij}}]$ {by the same arguments}.  The proof is completed. \hfill $\Box$

\subsection{A Multivariate Mat{\'e}rn Model}
\def\bomega{\boldsymbol{\omega}}
\def\bF{\boldsymbol{F}}
\begin{theorem} \label{Eugenio_y_Aninha_Forever}
	Let $\bF: [-1,1] \times [-1,1] \times [0,\infty) \to \R^{p \times p}$ be a continuous mapping such that \\
	1. $\bF(\cdot,\cdot, \omega)$ is positive definite on $\T^{d_1,d_2}$ for all $\omega \ge 0$; \\
	2. $\bF(0,0, \cdot)$ is pointwise integrable. \\
	Let 
	\begin{equation}
		\label{Eugenio} \bK(s,r,h):= \int_{\R^d} {\rm e}^{\mathsf{i } \langle \bomega, \bh \rangle } \bF (s,r, \|\bomega\|) {\rm d} \bomega.
	\end{equation}
	Then, $\bK $ is a positive definite matrix valued function on $\T^{d_1,d_2} \times \R^d$.
\end{theorem}
{\bf Proof.} By Assumption 1, and as a corollary of Theorem 3.4 of \cite{berg2016schoenberg}, we have that for every fixed $\omega$, the mapping $\bF$ can be uniquely written as 
$$ \bF (s,r,\omega) = 
\sum_{\bk=(k_1,k_2)\in \mathbb{N}_0^2} 
\bB_{\bk} (\omega) C_{k_1}^{(d_1-1)/2} (s) C_{k_2}^{(d_2-1)/2} (r), \qquad s,r \in [-1,1], \quad \omega \ge 0, $$
where, for every fixed $\omega$, $\{ \bB_{\bk}(\omega) \}_{\bk \in \mathbb{N}_0^2}$ is a sequence of positive definite matrices with the additional requirement that $\sum_{\bk} \bB_{\bk}(0)$ is finite. \\
We now rewrite the function $\bK$ as in (\ref{Eugenio}) as 
\begin{eqnarray*}
	\bK(s,r,h) &= &\int_{\R^d} {\rm e}^{\mathsf{i } \langle \bomega, \bh \rangle } \bF (s,r, \|\bomega\|) {\rm d} \bomega \\
	&=& \sum_{\bk} \int_{\R^d} {\rm e}^{\mathsf{i } \langle \bomega, \bh \rangle } \bB_{\bk}(\|\bomega\|) {\rm d} \bomega C_{k_1}^{(d_1-1)/2} (s) C_{k_2}^{(d_2-1)/2} (r) \\
	&=:& \sum_{\bk} \widetilde{\bB}_{\bk}(h) C_{k_1}^{(d_1-1)/2} (s) C_{k_2}^{(d_2-1)/2} (r)  \qquad s,r \in [-1,1], \quad h=\|\bh\|, \; \bh \in \R^d. 
\end{eqnarray*} 
By Cram{\'e}r's theorem \citep{cramer1940theory}, the functions $\widetilde{\bB}_{\bk}$ are positive definite in $\R^d$ for all $\bk \in \mathbb{N}_0^2$. Thus, we can invoke Theorem \ref{mauro_avoid!} and the proof will be completed if we show that the sequence $\{ \widetilde{\bB}_{\bk}(0) \}_{\bk}$ is convergent. This comes from Assumption 2 which guarantes that $\bK$ is bounded at $(0,0,0)$.
The proof is completed. \hfill $\Box$

For the following result, we need to define the Mat{\'e}rn family of functions $\mathcal{M}(\cdot; \alpha,\nu)$ through the identity 
\begin{equation}
	\label{Eugenio2} 
	\mathcal{M}(h; \alpha, \nu) = \frac{2^{1-\nu}}{\Gamma(\nu)} \left ( \alpha h \right )^{\nu} \mathcal{K}_{\nu} (\alpha h), \qquad h \ge 0, 
\end{equation}
for $\alpha, \nu$ being strictly positive, and where $\mathcal{K}_{\nu}$ denotes the modified Bessel function of the second kind of order $\nu>0$. Some works related to this family when analyzed over spheres can be found in \cite{terdik2015angular} and \cite{leonenko2017matern}.

\begin{theorem}
	Let $\big [ \nu_{ij}(\cdot,\cdot) \big ]_{i,j=1}^p$ be a matrix of functions $\nu_{ij}: [-1,1]^2 \to (0,\infty )$ that are continuous and such that $\nu_{ij}(\cdot,\cdot) = (\nu_{ii}(\cdot,\cdot)+\nu_{jj}(\cdot,\cdot))/2$. Let 
	\begin{equation}
		\label{Eugenio3}
		\rho_{ij}(s,r) = \beta_{ij} \frac{\Gamma(\nu_{ij}(s,r))}{\Gamma(\nu_{ij}(s,r)+d/2)}, \qquad s,r \in [-1,1], 
	\end{equation}
	where $[\beta_{ij}]$ is a matrix with $\beta_{ii}$=1. 
Assume that for any $0 < a \leq 1$,	the matrix-valued function $(s,r) \mapsto [ \beta_{ij} a^{\nu_{ij}(s,r)} ]$ is positive definite in $\T^{d_1,d_2}$.	
	Let $\bK: [-1,1] \times [-1,1] \times [0,\infty) \to \R^{p \times p}$ have elements $K_{ij}$ that are defined as 
	\begin{equation}
		\label{Eugenio4}  K_{ij}(s,r,h) = \sigma_i \sigma_j\rho_{ij}(s,r) \mathcal{M} \Big (h; \alpha, \nu_{ij}(s,r) \Big ),
	\end{equation}
with  $\sigma_i>0$, $i=1,\ldots,p$.	Then, $\bK$ is positive definite.   
\end{theorem}
{\bf Proof.} Let $$f(\omega; \alpha, \nu) =\frac{\Gamma(\nu+d/2) \alpha^{2 \nu} }{\Gamma(\nu) \pi^{d/2}} \left ( \alpha^2 + \omega^2 \right )^{-\nu-d/2}, \qquad \omega = \|\bomega\|, \; \bomega \in \R^d,  $$
with $\alpha$ and $\nu$ being strictly positive. The Fourier transform in $\R^d$ of the function $f$ is the function $\mathcal{M}(\cdot; \alpha, \nu)$, defined at (\ref{Eugenio2}). We now define 
$$ \widetilde{f}(\omega; \alpha, \nu) =  \frac{\Gamma(\nu)\pi^{d/2}}{\Gamma(\nu+d/2)} f (\omega; \alpha, \nu), \qquad \omega \ge 0, $$
having Fourier transform the function $$ \widetilde{\mathcal{M}}(h;\alpha,\nu) = \frac{\Gamma(\nu) \pi^{d/2}}{\Gamma(\nu+d/2)} \mathcal{M}(h; \alpha, \nu), \qquad h \ge 0. $$ 
 We consider the function $\bF:[-1,1] \times [-1,1] \times [0,\infty) \to \R^{p \times p}$ having elements $F_{ij}$ that are identically defined as 
\begin{eqnarray*}
	F_{ij} (s,r,\omega) &=& \sigma_i \sigma_ j 
	\beta_{ij} \widetilde{f} \left (\omega; {\alpha,} \nu_{ii}(s,r) \right )^{1/2} \widetilde{f} \left (\omega; {\alpha,} \nu_{jj}(s,r) \right )^{1/2}  \\
	&=& \sigma_i \sigma_ j \beta_{ij} \alpha^{2 \nu_{ij}(s,r)} \left ( \alpha^2 + \omega ^2 \right )^{-\nu_{ij}(s,r) - d/2 }. 
\end{eqnarray*} 
Clearly, $\bF $ satisfies the two hypothesis in Theorem \ref{Eugenio_y_Aninha_Forever}. 
 Thus, we can calculate the function $\bK$ as in (\ref{Eugenio}) to obtain a matrix valued mapping $\bK$ with elements $K_{ij}$ defined as in (\ref{Eugenio4}). The proof is completed by direct application of Theorem \ref{Eugenio_y_Aninha_Forever}.  \hfill $\Box$

\section{Evading from Stationarity and Isotropy in $\R^d$} \label{sec:4}

We have explored so far the constructions related to spherical isotropy on $\T^{d_1,d_2}$ coupled with Euclidean isotropy in $\R^d$. The next findings show characterizations related to relaxing this last hypothesis. Specifically, we refer to covariance functions $\bK_{\boldsymbol{Z}}$ such that 
\def\bk{\boldsymbol{k}}
\def\jj{j^{\prime}}
\def\buu{\bu^{\prime}}
\begin{equation}
	\label{menegatto2}  \bK_{\boldsymbol{Z}} \Big ( (\mathbf{x},\bu), (\mathbf{x}^{\prime},\bu^{\prime}) \Big)  = \left ( \Sigma^{1/2} \right )^{\top}  \bK( s,r, \bu,\bu^{\prime}  ) \Sigma^{1/2}, \qquad 
\end{equation} 
with $s = \langle x_1,{x^\prime_1} \rangle_1,  r = \langle x_2,{x^\prime_2} \rangle_2$, {$\mathbf{x}=(x_1,x_2),\mathbf{x}^\prime=(x_1^\prime,x_2^\prime)\in\T^{d_1,d_2}$,} $\bu,\bu^{\prime} \in \R^d$. The matrix $\Sigma^{1/2}$ is being defined around (\ref{menegatto}).  Under this setting, the mapping $\bK$ is defined over $[-1,1]^2 \times (\R^d)^2$ with image on $\R^{p \times p}$. Characterizations for this class have been elusive so far. We call $\Xi(\T^{d_1,d_2},\R^d)$ the class of mappings $\bK: \T^{d_1,d_2} \times \R^{d} \times  \T^{d_1,d_2} \times \R^{d} \to \R^{p \times p}$ such that
$$ \bK\Big ( (\bx, \bu),(\bx^{\prime},\bu^{\prime}) \Big ) = \sum_{\bk\in\mathbb{N}_0^2}
\sum_{j=1}^{\kappa_{d_1,k_1}} \sum_{\jj=1}^{\kappa_{d_2,k_2}} \Upsilon_{\bk,j,\jj}(\bu,\buu) Y_{k_1,j}(x_1)\overline{Y_{k_1,\jj}(x_1^{\prime})} Y_{k_2,j}(x_2)\overline{Y_{k_2,\jj}(x_2^{\prime})},  $$
for $s = \langle x_1,x_1^{\prime} \rangle_1$, $r=\langle x_2,x_2^{\prime} \rangle_2 $, $\bk=(k_1,k_2)$, $\{ Y_{k_i,j}\}_{k_i=0}^{\infty}$ denoting the set of spherical harmonics in $\S^{d_i}$, having a uniquely determined finite dimension $\kappa_{d_i,k_i}$ that  is 
	given by \eqref{apr1}.
 The expansion is uniformly convergent with respect to $(\bx,\bx^{\prime}) \in \T^{d_1,d_2}$. 
The next result is providing some insight into the knowledge of the class $\Xi(\T^{d_1,d_2},\R^d)$.
\begin{theorem} \label{thm:nonstat1}
	Let $\bK { :\T^{d_1,d_2} \times \R^{d} \times  \T^{d_1,d_2} \times \R^{d} \to}
	\R^{p \times p}$ belong to $\Xi(\T^{d_1,d_2},\R^d)$. 
	Then, \\
	(i) for every fixed $\bu,\bu^{\prime} \in \R^d$,
	$$ \sum_{\bk} \sum_j \sum_{\jj} \Big | \Upsilon_{\bk,j,\jj}^{(\ell,m)}(\bu,\buu) \Big |^2 < \infty, $$
	where $\Upsilon_{\bk,j,\jj}^{(\ell,m)}$ denotes the $(\ell,m)$-element of $\Upsilon_{\bk,j,\jj}$. Additionally, 
	$$ \Upsilon_{\bk,j,\jj}(\bu,\buu) = \int_{\T^{d_1,d_2}} \int_{\T^{d_1,d_2}} \bK\Big ( (\bx, \bu),(\bx^{\prime},\bu^{\prime}) \Big ) {\rm d} \lambda_{d_1,d_2} (\bx) {\rm d} \lambda_{d_1,d_2}(\bx^{\prime}), $$ 
	where ${\rm d} \lambda_{d_1,d_2} = {\rm d} \sigma_{d_1} \otimes {\rm d} \sigma_{d_2}$, with ${\rm d} \sigma_{d_i}$ denoting the Lebesgue measure on $\S^{d_i}$, $i=1,2$.\\
	(ii) $\bK$ is positive definite if and only if $\{ \Upsilon_{\bk, j, \jj } \}$ is a sequence of positive definite functions.
\end{theorem}

{\bf Proof.} Let $N,N^{\prime} \in \mathbb{N}$. Define 
\begin{eqnarray} 
	\label{eugenio_fazendo_pizza}
	&& \bK_{N,N^{\prime}}\Big ( (\bx, \bu),(\bx^{\prime},\bu^{\prime}) \Big )  \\ &=& \sum_{k_1=0}^{N} \sum_{k_2=0}^{N^{\prime}} \sum_{j=1}^{\kappa_{d_1,k_1}} \sum_{\jj=1}^{\kappa_{d_2,k_2}} \Upsilon_{\bk,j,\jj}(\bu,\buu) Y_{k_1,j}(x_1)\overline{Y_{k_1,\jj}(x_1^{\prime})} Y_{k_2,j}(x_2)\overline{Y_{k_2,\jj}(x_2^{\prime})}. \nonumber
\end{eqnarray}
Clearly, $\bK_{N,N^{\prime}}$ converges uniformly on $\T^{d_1,d_2} \times \T^{d_1,d_2}$ to $\bK$ as $N,N^{\prime}$ tend to infinity for any fixed $\bu,\bu^{\prime}$. This implies convergence in $\mathcal{L}_2(\T^{d_1,d_2} \times \T^{d_1,d_2}, \sigma_{d_1} \otimes \sigma_{d_2})$, the space of square integrable functions on $\T^{d_1,d_2} \times \T^{d_1,d_2}$ with respect to the tensor product measure $\sigma_{d_1} \otimes \sigma_{d_2}$. Here, square integrability is intended as finiteness of the Frobenius norm {$| \cdot |$}, so that we have
\begin{eqnarray*}
	\infty &>& \int_{\T^{d_1,d_2}} \int_{\T^{d_1,d_2}} \Bigg |\bK\Big ( (\bx, \bu),(\bx^{\prime},\bu^{\prime}) \Big )- \bK_{N,N^{\prime}}\Big ( (\bx, \bu),(\bx^{\prime},\bu^{\prime}) \Big )  \Bigg |^2 {\rm d} \sigma_{d_1}({\bx}) {\rm d} \sigma_{d_2}(\bx^{\prime}) \\
	&=&  \sum_{k_1=N}^{\infty} \sum_{k_2=N^{\prime}}^{\infty} \sum_{j{=1}}^{{\kappa_{d_1,k_1}}}
	 \sum_{\jj {=1}}^{{\kappa_{d_2,k_2}}} 
	 { \left|   \Upsilon_{\bk,j,\jj}(\bu,\buu) \right|^2},
\end{eqnarray*}
where the second line is due to legitimate inversion of series with integrals, in concert with the fact that the spherical harmonics are orthonormal basis in their respective spaces. Point (i) is established. Point (ii) can be proved by following the same arguments as in Theorem \ref{mauro_avoid!}. \hfill $\Box$

\section{Conclusions}
We have provided characterization theorems for the covariance functions of vector valued random fields that are continuously indexed over generalized spaces involving the hypertorus. In turn, we have proposed several parametric classes available in closed form. This work opens for several research directions. For instance, it might be extremely useful to study conditions for equivalence of Gaussian measures defined over compact subsets of these generalized spaces. This would have important consequences for maximum likelihood estimation, as well as for kriging prediction, under infill asymptotics \citep{yadrenko1983spectral, stein-book, bevilacqua2019estimation, arafat2018equivalence,bachoc2020asymptotically}. Such conditions would allow to evaluate the covariance functions proposed in this paper in terms of estimation and prediction.  Another relevant problem is related to modeling anisotropies in these generalized spaces, and how to estimate such structure.  \\
We should also mention that recent approaches allow to define reducibility over hyperspheres \citep{allard2016anisotropy, porcu2020reduction, senoussi2021nonstationary}. It would be definitely mandatory to study such approaches over the generalized spaces introduced in this paper.


\begin{thebibliography}{}
	
	\bibitem[Alegr\'{i}a et~al., 2021]{alegria2021mathcalffamily}
	Alegr\'{i}a, A., Cuevas-Pacheco, F., Diggle, P., and Porcu, E. (2021).
	\newblock The $\mathcal{F}$-family of covariance functions: A {Mat\'ern}
	analogue for modeling random fields on spheres.
	
	\bibitem[Allard et~al., 2016]{allard2016anisotropy}
	Allard, D., Senoussi, R., and Porcu, E. (2016).
	\newblock Anisotropy models for spatial data.
	\newblock {\em Mathematical Geosciences}, 48(3):305--328.
	
	\bibitem[Arafat et~al., 2018]{arafat2018equivalence}
	Arafat, A., Porcu, E., Bevilacqua, M., and Mateu, J. (2018).
	\newblock Equivalence and orthogonality of {Gaussian} measures on spheres.
	\newblock {\em Journal of Multivariate Analysis}, 167:306--318.
	
	\bibitem[Bachoc et~al., 2017]{bachoc2017gaussian}
	Bachoc, F., Gamboa, F., Loubes, J.-M., and Venet, N. (2017).
	\newblock A {Gaussian} process regression model for distribution inputs.
	\newblock {\em IEEE Transactions on Information Theory}, 64(10):6620--6637.
	
	\bibitem[Bachoc et~al., 2020]{bachoc2020asymptotically}
	Bachoc, F., Porcu, E., Bevilacqua, M., Furrer, R., and Faouzi, T. (2020).
	\newblock Asymptotically equivalent prediction in multivariate geostatistics.
	\newblock {\em arXiv preprint arXiv:2007.14684}.
	
	\bibitem[Berg et~al., 1984]{berg1984harmonic}
	Berg, C., Christensen, J. P.~R., and Ressel, P. (1984).
	\newblock {\em Harmonic analysis on semigroups: theory of positive definite and
		related functions}, volume 100.
	\newblock Springer.
	
	\bibitem[Berg et~al., 2018a]{berg1}
	Berg, C., Peron, A., and Porcu, E. (2018a).
	\newblock Orthogonal expansions related to compact {G}elfand pairs.
	\newblock {\em Exposithiones Mathematicae}, 36:259--277.
	
	\bibitem[Berg et~al., 2018b]{berg2}
	Berg, C., Peron, A., and Porcu, E. (2018b).
	\newblock Schoenberg's theorem for real and complex {H}ilbert spheres
	revisited.
	\newblock {\em Journal of Approximation Theory}, 228:58--78.
	
	\bibitem[Berg and Porcu, 2017]{berg2016schoenberg}
	Berg, C. and Porcu, E. (2017).
	\newblock From {S}choenberg coefficients to {S}choenberg functions.
	\newblock {\em Constructive Approximation}, 45(2):217--241.
	
	\bibitem[Bevilacqua et~al., 2019]{bevilacqua2019estimation}
	Bevilacqua, M., Faouzi, T., Furrer, R., and Porcu, E. (2019).
	\newblock Estimation and prediction using generalized {Wendland} covariance
	functions under fixed domain asymptotics.
	\newblock {\em Annals of Statistics}, 47(2):828--856.
	
	\bibitem[Cramer, 1940]{cramer1940theory}
	Cramer, H. (1940).
	\newblock On the theory of stationary random processes.
	\newblock {\em Annals of Mathematics}, 41(1):215--230.
	
	\bibitem[Estrade et~al., 2019]{estrade}
	Estrade, A., Fari{\~n}as, A., and Porcu, E. (2019).
	\newblock Covariance functions on spheres cross time: Beyond spatial isotropy
	and temporal stationarity.
	\newblock {\em Statistics $\&$ Probability Letters}, 151:1--7.
	
	\bibitem[Guella and Menegatto, 2016]{mene1}
	Guella, G. and Menegatto, V. (2016).
	\newblock Strictly positive definite kernels on a product of spheres.
	\newblock {\em Journal of Mathematical Analysis and Applications},
	435:286--301.
	
	\bibitem[Guella et~al., 2015]{mene2}
	Guella, J., Menegatto, V., and Peron, A. (2015).
	\newblock An extension of a theorem of {S}choenberg to products of spheres.
	\newblock {\em Banach Journal of Mathematical Analysis}, 435:286--301.
	
	\bibitem[Guella and Menegatto, 2020]{guella-men-cpd2020}
	Guella, J.~C. and Menegatto, V.~A. (2020).
	\newblock Conditionally positive definite matrix valued kernels on {E}uclidean
	spaces.
	\newblock {\em Constr. Approx.}, 52(1):65--92.
	
	\bibitem[Guella et~al., 2016]{mene4}
	Guella, J.~C., Menegatto, V.~A., and Peron, A.~P. (2016).
	\newblock Strictly positive definite kernels on a product of spheres {II}.
	\newblock {\em SIGMA Symmetry Integrability Geom. Methods Appl.}, 12:Paper No.
	103, 15.
	
	\bibitem[Guella et~al., 2017]{mene3}
	Guella, J.~C., Menegatto, V.~A., and Peron, A.~P. (2017).
	\newblock Strictly positive definite kernels on a product of circles.
	\newblock {\em Positivity}, 21(1):329--342.
	
	\bibitem[Leonenko and Malyarenko, 2017]{leonenko2017matern}
	Leonenko, N. and Malyarenko, A. (2017).
	\newblock Mat{\'e}rn class tensor-valued random fields and beyond.
	\newblock {\em Journal of Statistical Physics}, 168(6):1276--1301.
	
	\bibitem[Malyarenko and Ostoja-Starzewski, 2018]{malyarenko2018tensor}
	Malyarenko, A. and Ostoja-Starzewski, M. (2018).
	\newblock {\em Tensor-valued random fields for continuum physics}.
	\newblock Cambridge University Press.
	
	\bibitem[Marinucci and Peccati, 2011]{marinucci-peccati}
	Marinucci, D. and Peccati, G. (2011).
	\newblock {\em Random {F}ields on the {S}phere, {R}epresentation, {L}imit
		{T}heorems and {C}osmological {A}pplications}.
	\newblock Cambridge, New York.
	
	\bibitem[Mastrantonio et~al., 2016]{mastrantonio2}
	Mastrantonio, G., Jona~Lasinio, G., and Gelfand, A. (2016).
	\newblock Spatio-temporal circular models with non-separable covariance
	structure.
	\newblock {\em Test}, 25:331--350.
	
	\bibitem[Mastrantonio et~al., 2019]{mastrantonio}
	Mastrantonio, G., Jona~Lasinio, G., Pollice, A., Capotorti, G., Teodonio, L.,
	Genova, G., and Blasi, C. (2019).
	\newblock A hierarchical multivariate spatio-temporal model for clustered
	climate data with annual cycles.
	\newblock {\em Annals of Applied Statistics}, 13(2):797--823.
	
	\bibitem[Menegatto et~al., 2006]{2006P-valdir-claude-cpd-dotkernel}
	Menegatto, V.~A., Oliveira, C.~P., and Peron, A.~P. (2006).
	\newblock Conditionally positive definite dot product kernels.
	\newblock {\em J. Math. Anal. Appl.}, 321(1):223--241.
	
	\bibitem[Menegatto and Peron, 2004]{2004P-valdir-cpd-eucl}
	Menegatto, V.~A. and Peron, A.~P. (2004).
	\newblock Conditionally positive definite kernels on {E}uclidean domains.
	\newblock {\em J. Math. Anal. Appl.}, 294(1):345--359.
	
	\bibitem[Morimoto, 1998]{MR1641900}
	Morimoto, M. (1998).
	\newblock {\em Analytic functionals on the sphere}, volume 178 of {\em
		Translations of Mathematical Monographs}.
	\newblock American Mathematical Society, Providence, RI.
	
	\bibitem[Olver et~al., 2010]{olver2010nist}
	Olver, F.~W., Lozier, D.~W., Boisvert, R.~F., and Clark, C.~W. (2010).
	\newblock {\em NIST handbook of mathematical functions hardback and CD-ROM}.
	\newblock Cambridge university press.
	
	\bibitem[Porcu et~al., 2018]{PAF2016}
	Porcu, E., Alegr{\'i}a, A., and Furrer, R. (2018).
	\newblock Modeling spatially global and temporally evolving data.
	\newblock {\em International Statistical Review}, 86:344--377.
	
	\bibitem[Porcu et~al., 2016]{PBG16}
	Porcu, E., Bevilacqua, M., and Genton, M.~G. (2016).
	\newblock Spatio-temporal covariance and cross-covariance functions of the
	great circle distance on a sphere.
	\newblock {\em Journal of the American Statistical Association},
	111(514):888--898.
	
	\bibitem[Porcu et~al., 2019]{porcu2019axially}
	Porcu, E., Castruccio, S., Alegria, A., and Crippa, P. (2019).
	\newblock Axially symmetric models for global data: A journey between
	geostatistics and stochastic generators.
	\newblock {\em Environmetrics}, 30(1):e2555.
	
	\bibitem[Porcu et~al., 2020a]{porcu30}
	Porcu, E., Furrer, R., and Nychka, D. (2020a).
	\newblock 30 years of space--time covariance functions.
	\newblock {\em Wiley Interdisciplinary Reviews: Computational Statistics}, page
	e1512.
	
	\bibitem[Porcu et~al., 2020b]{porcu2020reduction}
	Porcu, E., Senoussi, R., Mendoza, E., and Bevilacqua, M. (2020b).
	\newblock Reduction problems and deformation approaches to nonstationary
	covariance functions over spheres.
	\newblock {\em Electronic Journal of Statistics}, 14(1):890--916.
	
	\bibitem[Porcu and White, 2020]{PW1}
	Porcu, E. and White, P. (2020).
	\newblock Random {F}ields on the {H}ypertorus. {P}art {I}: their {C}ovariance
	{M}odelling and {A}pplications.
	\newblock {\em Submitted}.
	
	\bibitem[Rudin, 1964]{rudin}
	Rudin, W. (1964).
	\newblock {\em Principles of mathematical analysis}, volume~3.
	\newblock McGraw-Hill New York.
	
	\bibitem[Schlather, 2010]{schlather2010some}
	Schlather, M. (2010).
	\newblock Some covariance models based on normal scale mixtures.
	\newblock {\em Bernoulli}, 16(3):780--797.
	
	\bibitem[Schoenberg, 1942]{schoen42}
	Schoenberg, I.~J. (1942).
	\newblock Positive definite functions on spheres.
	\newblock {\em Duke Math. J.}, 9(1):96--108.
	
	\bibitem[Senoussi and Porcu, 2021]{senoussi2021nonstationary}
	Senoussi, R. and Porcu, E. (2021).
	\newblock Nonstationary space--time covariance functions induced by dynamical
	systems.
	\newblock {\em Scandinavian Journal of Statistics}.
	
	\bibitem[Shirota and Gelfand, 2017]{shirota}
	Shirota, S. and Gelfand, A. (2017).
	\newblock {S}pace and circular time log {G}aussian {C}ox processes with
	application to crime event data.
	\newblock {\em Annals of Applied Statistics}, 11(2):481--503.
	
	\bibitem[Stein, 1999]{stein-book}
	Stein, M. (1999).
	\newblock {\em Interpolation of Spatial Data: Some Theory for Kriging}.
	\newblock Springer, New York.
	
	\bibitem[Szeg\H{o}, 1939]{szego}
	Szeg\H{o}, G. (1939).
	\newblock {\em Orthogonal {P}olynomials}, volume XXIII of {\em COLLOQUIUM
		PUBLICATIONS}.
	\newblock American Mathematical Society.
	
	\bibitem[Terdik, 2015]{terdik2015angular}
	Terdik, G. (2015).
	\newblock Angular spectra for non-{Gaussian} isotropic fields.
	\newblock {\em Brazilian Journal of Probability and Statistics},
	29(4):833--865.
	
	\bibitem[White and Porcu, 2019a]{white2019a}
	White, P. and Porcu, E. (2019a).
	\newblock Nonseparable covariance models on circles cross time: A study of
	{M}exico {C}ity ozone.
	\newblock {\em Environmetrics}, page e2558.
	
	\bibitem[White and Porcu, 2019b]{white}
	White, P. and Porcu, E. (2019b).
	\newblock Towards a complete picture of stationary covariance functions on
	spheres cross time.
	\newblock {\em Electronic Journal of Statistics}, 13:2566--2594.
	
	\bibitem[Yadrenko and Balakrishnan, 1983]{yadrenko1983spectral}
	Yadrenko, M.~I. and Balakrishnan, A.~V. (1983).
	\newblock {\em Spectral Theory of Random Fields (Spektral'naja Teorija
		Sluchajnykh Polej)}.
	\newblock Optimization Software, Publications Division.
	
\end{thebibliography}
\end{document}